\documentclass[preprint,review,12pt]{elsarticle}
\usepackage{amsmath}
\usepackage{amsfonts}
\usepackage{amssymb}
\usepackage{mathtools}
\usepackage[T1]{fontenc}
\usepackage{booktabs}
\usepackage{listings}
\usepackage{xcolor}
\usepackage{hyperref}
\usepackage{float}
\usepackage{graphicx}
\usepackage{amsthm}
\usepackage{lineno}
\theoremstyle{plain}
\newtheorem{theorem}{Theorem}
\newtheorem{corollary}{Corollary}
\theoremstyle{definition}
\newtheorem{definition}{Definition}
\newcommand{\Hupper}{H_{\mathrm{upper}}}
\newcommand{\Hlower}{H_{\mathrm{lower}}}
\newcommand{\Hin}{H_{\mathrm{in}}}
\newcommand{\Hout}{H_{\mathrm{out}}}
\hypersetup{colorlinks=true, linkcolor=blue, urlcolor=cyan}
\journal{Information Sciences}
\begin{document}
\begin{frontmatter}
\title{Compact Slack-Free Penalty Functions for Inequality
Constraints in Higher-Order Unconstrained Binary Optimization}
\author[sm,qf]{Naoaki Mochida}
\author[ist,vc,qf]{Shoya Yasuda\corref{cor1}}
\ead{yasuda@vigne-cla.com}
\cortext[cor1]{Corresponding author.}

\affiliation[sm]{organization={St.\ Mary's International School},
  addressline={1-6-19 Seta, Setagaya-ku},
  city={Tokyo}, postcode={158-8668}, country={Japan}}
\affiliation[qf]{organization={Quantum Fabric Inc.},
  addressline={2-24-12 Shibuya, Shibuya-ku},
  city={Tokyo}, postcode={150-6139}, country={Japan}}
\affiliation[ist]{organization={Institute of Science Tokyo, School of Environment and Society},
  addressline={2-12-1 Ookayama, Meguro-ku},
  city={Tokyo}, postcode={152-8550}, country={Japan}}
\affiliation[vc]{organization={Vignette \& Clarity, Inc.},
  addressline={5-21 Shirahataminami-cho, Kanagawa-ku},
  city={Yokohama}, postcode={221-0073}, country={Japan}}

\begin{abstract}
Inequality constraints in higher-order unconstrained binary optimization (HUBO) are commonly handled using slack variables or by enumerating infeasible assignments, which respectively enlarge the search space or produce penalties whose size grows with the number of forbidden values. We construct exact slack-free penalty functions for upper- and lower-bound constraints on binary-encoded integers. Each function is constructed hierarchically from the threshold bits through a recursive addition-and-multiplication rule, and its correctness for arbitrary bit widths and representable thresholds is proved by mathematical induction. Equivalent explicit finite sum-product representations are also derived, enabling direct analysis of product-term count and polynomial degree. Interval constraints are obtained by adding the lower- and upper-bound penalties, while exclusion constraints are obtained by multiplying them. Compared with the trivial enumerative construction, the proposed method reduces the worst-case product-term count from linear to logarithmic in the representable range. Except at one non-vacuous boundary threshold for each bound, it also reduces both product-term count and degree. The resulting single-bound penalties have minimum positive value 1
and maximum value at most the bit width, which may simplify penalty-weight selection. The same construction also extends exactly
to fixed-scale fractional variables through appropriate scaling of
the integer thresholds.
\end{abstract}


\end{frontmatter}

\section*{Glossary}

\begin{description}
\item[HUBO] Higher-Order Unconstrained Binary Optimization.
\item[QUBO] Quadratic Unconstrained Binary Optimization.
\item[Penalty function] A nonnegative function that is zero on feasible assignments and strictly positive on infeasible assignments.
\item[Literal] A binary variable $x_i$ or its complement $1-x_i$.
\item[Product term] A product of literals.
\item[Product-term count] The number of product terms in the generated finite sum-product representation.
\item[Degree] The maximum number of literals contained in a product term.
\end{description}

\section{Introduction}
\label{sec:intro}
\subsection{Motivation: Inequality Constraints in Quantum Optimization}
\label{subsec:motivation}

Combinatorial optimization is one of the most promising near-term applications of quantum computing, including quantum annealing~\cite{Kadowaki1998,Hauke2020,Yarkoni2022} and the Quantum Approximate Optimization Algorithm (QAOA)~\cite{Farhi2014,Abbas2024,Blekos2024,Hadfield2019}. Such problems can also be formulated as Higher-Order Unconstrained Binary Optimization (HUBO), also referred to as Polynomial Unconstrained Binary Optimization (PUBO) or Higher-Order Binary Optimization (HOBO), which generalizes Quadratic Unconstrained Binary Optimization (QUBO) by allowing interactions of degree higher than two~\cite{Abbas2024,Lucas2014}. Within this framework, equality and inequality constraints are not on equal footing. An equality constraint of the form $f=0$ can be incorporated as a penalty term $f^{2}$, without introducing any variable beyond those already present in $f$~\cite{Glover2019}. However, in inequality constraints, the standard treatments instead rely on slack variables or one-hot encodings, both of which enlarge the search space with auxiliary variables that do not belong to the original problem. The reliance on auxiliary variables makes inequality constraints one of the central difficulties in conventional QUBO and HUBO formulations.

QUBO formulations already make use of binary encoding---representing an integer through its positional base-2 expansion, where each bit contributes a successive power of two, rather than allocating one variable per possible value---to represent a bounded integer variable compactly, requiring only $O(\log R)$ qubits for a variable with bound $R$~\cite{KarimiRonagh2019,Tamura2021}, rather than $O(R)$ qubits as in a one-hot representation. Binary encoding, however, only fixes how an integer is represented as bits. It supplies no mechanism, by itself, for penalizing a violation of an inequality on that integer. This gap is independent of whether the underlying model is QUBO or HUBO. Regardless of the degree of interactions allowed, an explicit penalty construction is still required to make every infeasible bit-string energetically unfavorable~\cite{Glover2019,Abbas2024}. Such overhead is especially consequential on current noisy intermediate-scale quantum (NISQ) hardware~\cite{Preskill2018,Abbas2024}, where any additional structure imposed on the penalty function may affect the qubit count and connectivity needed for implementation.

A comparable intuitive approach exists in HUBO, which we call as 
\emph{trivial penalty function}, which enumerates the infeasible
assignments one by one and therefore does not exploit the higher-order
setting to reduce the underlying enumerative structure. Although this
construction is valid, it has several structural disadvantages. To see
why a \emph{trivial} penalty function is problematic, consider the most
direct way of enforcing an upper-bound constraint $X \leq U$, where
$X$ denotes the integer determined by $k$ binary digits through the
standard binary expansion introduced in
Section~\ref{subsec:prelim} (Eq.~\eqref{eq:binary_expansion}), so that
$X$ ranges over $\{0,1,\ldots,2^k-1\}$ rather than being an independent
free parameter. The formulation enumerates every forbidden integer
$X>U$ and adds one term per forbidden bit string. For example, for
$k=3$ and $U=4$, the integers $5=(101)_2$, $6=(110)_2$, and
$7=(111)_2$ are forbidden, and the \emph{trivial} penalty function
\begin{equation}
  H = x_3(1-x_2)x_1 + x_3x_2(1-x_1) + x_3x_2x_1
  \label{eq:trivial}
\end{equation}
correctly penalizes exactly these three configurations. The construction
is, in this sense, correct. Yet, it carries at least four compounding
drawbacks.

\textbf{(i) Growth in product-term count.} \mbox{}\\
The trivial construction adds exactly one product term for each forbidden
integer, so its product-term count equals the number of integers that
violate the threshold, namely $2^k-U-1$. For fixed $k$, this count grows
as $U$ decreases and, in the worst case, scales linearly with the number
of values representable by $X$. For classical or simulated higher-order
optimization, every additional product term increases the work required
to evaluate the penalty function or its local energy changes. On quantum
hardware, a larger collection of higher-order product terms may also
increase the resources required for quadratization, embedding, or
physical implementation, depending on the available interaction
structure~\cite{Abbas2024,Konz2021}.

\textbf{(ii) Fixed maximum degree.} \mbox{}\\
Each enumerative product term contains one literal for each of the
$k$ binary variables and therefore has degree $k$. Consequently,
whenever the constraint excludes at least one bit string, the maximum
product-term degree of the trivial penalty function is exactly $k$.
In classical and simulated higher-order optimization, such terms require
the evaluation of interactions involving all $k$ literals, although the
actual computational cost depends on the implementation and the use of
incremental updates. On current quantum annealers based on two-body
interactions, such higher-order terms must generally be reduced to
quadratic form. This reduction may introduce auxiliary variables and
additional embedding overhead, depending on the chosen quadratization
method and the hardware
connectivity~\cite{Mandal2020,Dattani2019,Konz2021,Abbas2024}.
Higher-order terms can therefore impose a greater computational or
implementation burden in both simulated and quantum settings.

\textbf{(iii) No use of hierarchical structure.}\mbox{}\\
The proposed upper bound and lower bound constructions obtain the
penalty function for $k+1$ bits from the corresponding function for
$k$ bits by adding one further hierarchical layer. The \emph{trivial}
penalty function instead treats every infeasible assignment
independently and shares no substructure among the resulting product
terms. This lack of reuse is a structural limitation because it
prevents a larger instance from being constructed from components
already obtained for fewer bits. It also contrasts with the
hierarchical organization commonly used in quantum computing, where
more complex operations are constructed from simpler operations, as
illustrated by decompositions of multiqubit controlled gates into
elementary gates~\cite{Barenco1995}. The proposed recursion preserves
an analogous hierarchy at the level of the penalty function before
either an annealing implementation or a gate model implementation is
chosen.

\textbf{(iv) Poor mathematical reusability and extensibility.}\mbox{}\\
This last drawback is not a matter of computational cost but of mathematical structure. A general justification of the trivial construction is not, in itself, hard to obtain. For any fixed $k$ and $U$, one can verify that the enumeration in Eq.~\eqref{eq:trivial} contributes exactly one term per forbidden integer and penalizes no feasible one. What the construction fails to do is convert the binary structure of the threshold $U$ into a compact recursive rule. Since every forbidden pattern is enumerated independently, there is no recursive relation carrying the penalty construction at bit width $k$ over to bit width $k+1$, and no feature of $U$'s own binary expansion is exploited beyond its role in listing the forbidden integers. This limitation becomes more consequential when the upper-bound penalty is combined with a lower-bound penalty to form an interval constraint, or with its complement to form an exclusion constraint. The component penalties still originate from term-by-term enumeration, so their composition inherits the same enumeration-sized structure rather than becoming more compact. The trivial construction is therefore less suitable for systematic mathematical reuse and extension, even though it is not, on its own terms, incorrect.

Addressing these four limitations therefore calls for a non-trivial penalty construction that exploits the binary structure of the threshold rather than enumerating infeasible assignments. Such a construction should reduce the generated product-term count from $O(R)$ in the representable range size $R=2^k$ to $O(\log R)$, suppress the maximum degree whenever possible, reuse the penalty structure already obtained for lower bit widths, and support systematic extension to compound constraints through arithmetic operations. A penalty construction with these properties would provide a more compact, hierarchical, and mathematically extensible alternative to the trivial enumerative formulation.
In this paper, we present an explicit hierarchical slack-free construction for inequality constraints in higher-order unconstrained binary optimization, implemented through recursive rules on the threshold bits.

\subsection{Contributions}
\label{subsec:contributions}

We make three principal contributions.

First, for any $k$-bit integer variable $X$ and any threshold
$U$ satisfying $0 \leq U < 2^k$, we introduce a hierarchical construction of an exact upper-bound penalty function
$\Hupper(X,U)$. We prove by mathematical induction that the
polynomial generated by this construction vanishes if and only if
$X \leq U$ for arbitrary bit widths and thresholds. We also derive
an equivalent explicit finite sum-product representation from the
recursion. No slack variables are introduced.

Second, for any $k$-bit integer variable $X$ and any threshold
$L$ satisfying $0 \leq L < 2^k$, we introduce a corresponding
hierarchical construction of an exact lower-bound penalty function
$\Hlower(X,L)$. We prove by mathematical induction that the
polynomial generated by this construction vanishes if and only if
$X \geq L$ for arbitrary bit widths and thresholds. We also derive
an equivalent explicit finite sum-product representation and
establish a one's-complement duality with the upper-bound
construction. This duality acts by complementing the encoded
integer and the threshold, rather than by applying a Boolean
operation to the penalty function.

Third, we use the non-negativity of the upper-bound and lower-bound
penalties to derive exact composite constraints. An interval
constraint $L \leq X \leq U$ is represented by
\[
\Hlower(X,L)+\Hupper(X,U),
\]
whereas an exclusion constraint $X \leq L$ or $U \leq X$ is represented
by
\[
\Hupper(X,L)\Hlower(X,U).
\]

Together, these results replace enumeration-sized inequality penalties with hierarchical constructions whose product-term count is bounded by the bit width, while also reducing the polynomial degree over nearly the entire non-vacuous threshold range. They also provide a unified arithmetic framework for constructing single-bound, interval, and exclusion penalties directly from the binary structure of the thresholds.

\subsection{Related Work and Positioning}
\label{subsec:related}

Standard QUBO formulations commonly encode inequalities through slack variables, which enlarge the binary search space~\cite{Glover2019,Vyskocil2019,Abbas2024}.
This approach enlarges the search space because the slack variable
must also be encoded by binary variables~\cite{Tamura2021,CeselliPremoli2023,Bontekoe2023}.

Several methods reduce or avoid this auxiliary-variable overhead.
Unbalanced penalization uses an asymmetric penalty in place of the
standard slack-variable construction~\cite{MontanezBarrera2023,MontanezBarrera2024}. Augmented-Lagrangian
and related multiplier methods instead enforce constraints through
iterative updates of multipliers or penalty parameters~\cite{Karimi2017,Djidjev2023,Gabbassov2025}. These approaches operate
at the algorithmic or iterative-reformulation level. Our objective is
different. We instead construct an exact polynomial penalty that is determined directly from a given binary threshold and requires no iterative parameter updates.

Slack-free penalties have also been studied in gate-based quantum optimization, including exponential and step-like constructions that suppress infeasible states without explicitly encoding slack variables~\cite{LeeAndLau2025}. In contrast, the present work focuses on exact polynomial penalty functions defined directly on the binary variables of threshold constraints. More recent slack-free approaches have continued in this direction.
Lee and Lau~\cite{LeeLau2026} combine nonlinear custom penalties with
finite-sampling estimation and CVaR in gate-based variational
optimization, whereas Hagiwara et al.~\cite{Hagiwara2026} combine
unbalanced penalization with an iterative Ohzeki-type update and deep
unfolding. These approaches address slack elimination at the
variational-sampling or iterative-algorithm level. In contrast, the
present work concerns the algebraic construction of an exact
finite pseudo-Boolean penalty polynomial from the bits of a fixed
threshold, without iterative updates in the penalty construction.

Two lines of prior work are particularly relevant, albeit from different perspectives. In the higher-order binary-optimization setting, the work of Minato~\cite{Minato2025} is the closest in problem setting because it directly addresses inequality constraints for binary-encoded integers in a HOBO formulation. Although this preprint was subsequently withdrawn, with the author noting overly restrictive assumptions and a critical error, we discuss it because of its close relation to the present problem setting. Its construction assigns separate product terms to forbidden bit patterns and therefore corresponds to the trivial enumerative penalty used as the baseline in this paper. In classical integer optimization, a related mathematical structure appears in the study of lexicographically bounded binary vectors. Muldoon et al.~\cite{Muldoon2013} used minimal-cover and set-covering inequalities to characterize binary vectors satisfying prescribed lexicographic upper or lower bounds and applied these results to base-2 representations of bounded integer variables. Related convex-hull characterizations were subsequently developed by Adams et al.~\cite{AdamsBelottiShen2016}. These works exploit the binary structure of a threshold to characterize the feasible set through systems of valid inequalities or polyhedral descriptions. The novelty of the present work is therefore not the characterization of the lexicographically bounded binary set itself, but the solution of a different mathematical problem required by unconstrained higher-order binary optimization: constructing, directly from the threshold bits, an exact nonnegative pseudo-Boolean penalty polynomial whose zero-energy set coincides exactly with the feasible set and whose value is strictly positive for every infeasible assignment. Unlike an external system of inequalities, such a penalty can be incorporated directly into an unconstrained higher-order objective as a constraint-energy term. We further provide hierarchical upper- and lower-bound constructions based on recursive threshold-bit rules that require no auxiliary variables in their original higher-order form, derive explicit finite sum-product representations that permit direct characterization of product-term count and degree, and exploit nonnegativity to compose the resulting penalties arithmetically into exact interval and exclusion constraints. Thus, the contribution is not merely a reformulation of a known binary threshold condition, but a hierarchical and composable zero-energy polynomial realization of that condition for higher-order unconstrained binary optimization.

Within the higher-order penalty-function setting considered here,
the structural difference from the enumerative baseline is
substantial. Let $R=2^k$ denote the size of the representable
integer range. The trivial penalty construction contains one
product term for each forbidden integer and can therefore require
$O(R)$ product terms. By contrast, each proposed
single-bound penalty is generated directly from the $k$ threshold
bits and contains at most $k=O(\log R)$ generated product terms.
The resulting product-term and degree comparisons are established in
Section~\ref{subsec:degree_count}.

Threshold constraints on binary-encoded integers have also been
studied in Boolean satisfiability~\cite{Warners1998,EenSorensson2006,Abio2012}. SAT represents
feasibility through Boolean formulas or clauses, whereas the present
approach represents a threshold constraint by a nonnegative
pseudo-Boolean penalty function whose zero set is exactly the feasible
set. In quantum optimization, such a function may also be interpreted
as a constraint Hamiltonian~\cite{Boros2002}. This difference in
representation leads naturally to the arithmetic composition rules
developed in Section~\ref{sec:logic}.

\section{Main Results: Construction and Proof}
\label{sec:main}

\subsection{Preliminaries and Notation}
\label{subsec:prelim}

A $k$-bit binary-encoded integer variable $X$ represents an integer in
$\{0,1,\ldots,2^k-1\}$ through
\begin{equation}
    X=\sum_{i=1}^{k}x_i2^{i-1},
    \qquad x_i\in\{0,1\},
    \label{eq:binary_expansion}
\end{equation}
and is written as
\[
    X=(x_kx_{k-1}\cdots x_2x_1)_2.
\]
A \emph{penalty function} is a nonnegative function
that is zero for every feasible binary assignment and strictly
positive for every infeasible assignment. Every penalty function
considered in this paper is expressed as a polynomial in the original
binary variables~\cite{Boros2002}.

We use the term penalty function because the constructions developed here are mathematical functions independent of a particular optimization implementation. In quantum optimization, they may be incorporated as constraint terms in a Hamiltonian. The same functions can also be used in quantum-inspired and general higher-order binary optimization solvers. We nevertheless retain the
notation $H$, including $\Hupper$, $\Hlower$, $\Hin$, and $\Hout$,
throughout the paper.

We call a penalty construction \emph{trivial} when it enumerates
the infeasible binary assignments and introduces a separate product
term for each of them. In contrast, we call a construction
\emph{non-trivial} when it generates the penalty directly from the
binary digits of the threshold without enumerating the infeasible
assignments. The constructions proposed in this paper are
non-trivial in this sense.

\begin{definition}[The four canonical penalty functions]
\label{def:four}
Let $X$ be the integer encoded by Eq.~\eqref{eq:binary_expansion}, and
let $L,U\in\{0,1,\ldots,2^k-1\}$. The four penalty functions considered
in this paper satisfy
\begin{align}
    \Hupper(X,U)=0
    &\quad\Longleftrightarrow\quad
    X\leq U,
    \label{eq:def-upper}\\
    \Hlower(X,L)=0
    &\quad\Longleftrightarrow\quad
    L\leq X,
    \label{eq:def-lower}\\
    \Hin(X,L,U)=0
    &\quad\Longleftrightarrow\quad
    L\leq X\leq U,
    \qquad L\leq U,
    \label{eq:def-in}\\
    \Hout(X,L,U)=0
    &\quad\Longleftrightarrow\quad
    X\leq L\ \text{or}\ U\leq X,
    \qquad L\leq U.
    \label{eq:def-out}
\end{align}
\end{definition}

Definition~\ref{def:four} specifies the required zero set of each
penalty function. The symbols $\Hupper$, $\Hlower$, $\Hin$, and
$\Hout$ denote the particular constructions developed in the
following subsections.

For a $k$-bit integer, the upper-bound constraint $X\leq U$ is
\emph{vacuous} when $U=2^k-1$, because every representable value of
$X$ is feasible. Similarly, the lower-bound constraint $L\leq X$ is
vacuous when $L=0$. The corresponding penalty function is the zero
polynomial. The zero polynomial has product-term count zero and no
assigned degree. Accordingly, vacuous endpoint constraints are
included in product-term counts with value zero but are excluded
from degree comparisons.

A \emph{literal} is either $x_i$ or $1-x_i$. The latter is a complemented literal of the same binary variable and does not introduce an additional variable.

A \emph{product term} is a product of literals. A penalty function is
regarded as a finite sum of such product terms, and its
\emph{product-term count} is the number of products in that sum.
Common factors written outside parentheses do not change this count.
For example,
\begin{equation}
    x_6x_5(x_4+x_3x_2)
    =
    x_6x_5x_4+x_6x_5x_3x_2
\end{equation}
contains two product terms. In contrast,
\begin{equation}
    (1-x_2)(1-x_1)
\end{equation}
contains one product term because it is a single product of two
literals. A complemented literal is not expanded into a constant and
a variable when the product-term count is evaluated. This convention is applied equally to both constructions and counts the product terms in their generated sum-product representations, not terms created only by algebraic expansion.

The \emph{degree} of a product term is the number of literals contained
in that product. The degree of a nonzero penalty function is the
maximum degree among its product terms. For example,
$x_6x_5x_4$ and $(1-x_6)(1-x_5)(1-x_4)$ both have degree three. The
function $x_6x_5(x_4+x_3x_2)$ has degree four because its two product
terms have degrees three and four. Because every literal is affine in a single binary variable, the literal-product degree used here is consistent with the ordinary total degree of each generated multilinear product term.

Throughout this paper, product-term count and degree are evaluated on
the expressions generated by the proposed construction and
the trivial construction. Both quantities are measured before
algebraic simplification, collection of identical terms, expansion of
complemented literals, or quadratization.

\subsection{Upper-Bound Penalty Function
  \texorpdfstring{$\Hupper(X,U)$}{H\textunderscore upper(X,U)}}
\label{subsec:upper}

We now construct the upper-bound penalty function $\Hupper(X,U)$ for $0 \leq U < 2^k$.

\begin{theorem}[Exact recursive upper-bound penalty]
\label{thm:upper}
For every integer $k\geq 1$ and every threshold
$U\in\{0,1,\ldots,2^k-1\}$, define the $k$-bit upper-bound penalty $H_k(X,U) \equiv H_{\mathrm{upper}}(X,U)$ recursively as
follows. For $k=1$,
\[
    H_1(X,0)=x_1,
    \qquad
    H_1(X,1)=0.
\]
For the extension from $k$ to $k+1$ bits, write
\[
    X_{k+1}=2^k x_{k+1}+X_k
\]
and
\[
    U=2^k u_{k+1}+U',
\]
where $x_{k+1},u_{k+1}\in\{0,1\}$ and
$0\leq X_k,U'<2^k$, and define
\[
H_{k+1}(X_{k+1},U)
=
\begin{cases}
x_{k+1}+H_k(X_k,U'),
& u_{k+1}=0,\\[4pt]
x_{k+1}H_k(X_k,U'),
& u_{k+1}=1.
\end{cases}
\]
Then $H_k(X,U)$ is a nonnegative integer-valued polynomial in the
original binary variables satisfying
\[
    H_k(X,U)=0
    \quad\Longleftrightarrow\quad
    X\leq U.
\]

\end{theorem}

\begin{proof}
We prove the stated properties by mathematical induction on the
bit width $k$.

\textit{Base case.}
Let $k=1$. Then $X=x_1$, and the possible thresholds are $U=0$
and $U=1$.

For $U=0$, the recursive definition gives
\[
    H_1(X,0)=x_1.
\]
Hence $H_1(X,0)=0$ exactly when $X=0$, and it is strictly positive
exactly when $X=1$.

For $U=1$, the recursive definition gives
\[
    H_1(X,1)=0.
\]
Both representable values satisfy $X\leq 1$, so the required
property also holds for $U=1$. Therefore, the theorem holds for
$k=1$.

\textit{Inductive hypothesis.}
Assume that for some $k\geq 1$ and for every threshold
$U'\in\{0,1,\ldots,2^k-1\}$, the recursively constructed
polynomial $H_k(X_k,U')$ is nonnegative and integer-valued and
satisfies
\[
    H_k(X_k,U')=0
    \quad\Longleftrightarrow\quad
    X_k\leq U',
\]
and
\[
    H_k(X_k,U')>0
    \quad\Longleftrightarrow\quad
    X_k>U'.
\]

\textit{Inductive step.}
Consider an arbitrary $(k+1)$-bit integer and an arbitrary
$(k+1)$-bit threshold written as
\[
    X_{k+1}=2^kx_{k+1}+X_k
\]
and
\[
    U=2^ku_{k+1}+U',
\]
where $x_{k+1},u_{k+1}\in\{0,1\}$ and
$0\leq X_k,U'<2^k$.

\textit{Case A: $u_{k+1}=0$.}
Then $U=U'$. By the recursive definition,
\begin{equation}
    H_{k+1}(X_{k+1},U)
    =
    x_{k+1}+H_k(X_k,U').
    \label{eq:upper_caseA}
\end{equation}
Both terms are nonnegative. Therefore,
\[
    H_{k+1}(X_{k+1},U)=0
\]
if and only if
\[
    x_{k+1}=0
    \quad\text{and}\quad
    H_k(X_k,U')=0.
\]
By the inductive hypothesis, this is equivalent to
\[
    x_{k+1}=0
    \quad\text{and}\quad
    X_k\leq U'.
\]
Since $U=U'<2^k$, these conditions are equivalent to
\[
    X_{k+1}\leq U.
\]
If either condition fails, at least one of the two nonnegative
terms is strictly positive, and hence
\[
    H_{k+1}(X_{k+1},U)>0.
\]
Thus,
\[
    H_{k+1}(X_{k+1},U)=0
    \quad\Longleftrightarrow\quad
    X_{k+1}\leq U
\]
in Case A.

\textit{Case B: $u_{k+1}=1$.}
Then $U=2^k+U'$. By the recursive definition,
\begin{equation}
    H_{k+1}(X_{k+1},U)
    =
    x_{k+1}H_k(X_k,U').
    \label{eq:upper_caseB}
\end{equation}

If $x_{k+1}=0$, then
\[
    X_{k+1}=X_k<2^k\leq U,
\]
and simultaneously
\[
    H_{k+1}(X_{k+1},U)=0.
\]

If $x_{k+1}=1$, then
\[
    H_{k+1}(X_{k+1},U)
    =
    H_k(X_k,U'),
\]
while
\[
    X_{k+1}\leq U
    \quad\Longleftrightarrow\quad
    X_k\leq U'.
\]
By the inductive hypothesis,
\[
    H_{k+1}(X_{k+1},U)=0
    \quad\Longleftrightarrow\quad
    X_{k+1}\leq U.
\]
Likewise,
\[
    H_{k+1}(X_{k+1},U)>0
    \quad\Longleftrightarrow\quad
    X_{k+1}>U.
\]
Thus the required property holds in Case B.

In both cases, $H_{k+1}$ is obtained from the nonnegative
integer-valued polynomial $H_k$ by addition or multiplication with
the binary variable $x_{k+1}$. Hence it remains a nonnegative
integer-valued polynomial, and its zero set is exactly
$\{X_{k+1}:X_{k+1}\leq U\}$.

Therefore, the stated properties hold for bit width $k+1$ whenever
they hold for bit width $k$. By mathematical induction, the theorem
holds for every integer $k\geq 1$.
\end{proof}

\begin{corollary}[Explicit finite sum-product representation]
\label{cor:upper_explicit}
Let
\[
    U=(u_ku_{k-1}\cdots u_1)_2.
\]
Expanding the recurrence in Theorem~\ref{thm:upper} gives
\begin{equation}
    \Hupper(X,U)
    =
    \sum_{\substack{1\leq i\leq k\\u_i=0}}
    x_i
    \prod_{\substack{j>i\\u_j=1}}x_j.
    \label{eq:upper_explicit}
\end{equation}
An empty product is defined as one, and an empty sum is defined as
zero. Each zero bit $u_i=0$ of the threshold therefore corresponds to
exactly one product term.
\end{corollary}
\begin{proof}
Starting from the one-bit construction and successively appending higher-order threshold bits, a zero threshold bit creates a new product term initialized by the corresponding variable, whereas a one threshold bit multiplies every existing product term by the corresponding variable. Repeated application of these two rules gives
Eq.~\eqref{eq:upper_explicit}.
\end{proof}

As a concrete illustration of the base case and one inductive step,
Table~\ref{tab:upper_k2} lists the resulting penalty functions for all
two-bit upper thresholds.

\begin{table}[H]
\centering
\caption{Upper-bound penalty functions for all thresholds at $k=2$.}
\label{tab:upper_k2}
\begin{tabular}{cc}
\toprule
Upper threshold $U$ & $\Hupper(X,U)$ \\
\midrule
$0=(00)_2$ & $x_2+x_1$ \\
$1=(01)_2$ & $x_2$ \\
$2=(10)_2$ & $x_2x_1$ \\
$3=(11)_2$ & $0$ \\
\bottomrule
\end{tabular}
\end{table}

The table illustrates the zero-set property established in
Theorem~\ref{thm:upper}.

Table~\ref{tab:upper_k3} extends the construction to three bits. Each
entry is obtained from the corresponding two-bit function by
applying Case A or Case B of the inductive step.

\begin{table}[H]
\centering
\caption{Upper-bound penalty functions for all thresholds at $k=3$.}
\label{tab:upper_k3}
\begin{tabular}{ccc}
\toprule
$(u_3u_2u_1)_2$ & $U$ & $\Hupper(X,U)$ \\
\midrule
000 & 0 & $x_3+x_2+x_1$ \\
001 & 1 & $x_3+x_2$ \\
010 & 2 & $x_3+x_2x_1$ \\
011 & 3 & $x_3$ \\
100 & 4 & $x_3(x_2+x_1)$ \\
101 & 5 & $x_3x_2$ \\
110 & 6 & $x_3x_2x_1$ \\
111 & 7 & $0$ \\
\bottomrule
\end{tabular}
\end{table}

Each function in Table~\ref{tab:upper_k3} has the required zero set.
The table also shows that the construction reuses the penalty function
formed for the lower bit positions. At the boundary threshold
$U=6=2^3-2$, the proposed function is $x_3x_2x_1$ and coincides with
the trivial construction.

\paragraph*{Algorithm~1: Iterative construction of
  \texorpdfstring{$\Hupper$}{H\_upper}}\mbox{}\\
Given the binary representation
\[
    U=(u_ku_{k-1}\cdots u_2u_1)_2,
\]
the upper-bound penalty function $\Hupper(X,U)$ is constructed as
follows.
\begin{enumerate}
    \item If $u_1=0$, initialize $H = x_1$.
          If $u_1=1$, initialize $H = 0$.

    \item For $i=2,3,\ldots,k$, update
    \[
        H\leftarrow
        \begin{cases}
            x_i+H, & \text{if } u_i=0,\\
            x_iH,  & \text{if } u_i=1.
        \end{cases}
    \]
\end{enumerate}
The final value of $H$ is $\Hupper(X,U)$.

Algorithm 1 is the iterative form of the recurrence in Theorem~\ref{thm:upper}, and its correctness therefore follows immediately from the theorem.

\paragraph*{Example $(k=6,\ U=53)$} \mbox{}\\
For
\[
    U=53=(110101)_2,
\]
Algorithm~1 gives
\begin{equation}
    \Hupper(X,53)
    =
    x_6x_5(x_4+x_3x_2).
    \label{eq:example_upper}
\end{equation}
This function has product-term count two and degree four. In contrast,
the trivial construction contains
\[
    2^6-53-1=10
\]
product terms and has degree six. Thus, the proposed construction
reduces the product-term count from ten to two and the degree from six
to four.

\subsection{Lower-Bound Penalty Function
  \texorpdfstring{$\Hlower(X,L)$}{H\textunderscore lower(X,L)}}
\label{subsec:lower}

We now construct the lower-bound penalty function
$\Hlower(X,L)$ for $0\leq L<2^k$.

\begin{theorem}[Exact recursive lower-bound penalty]
\label{thm:lower}
For every integer $k\geq 1$ and every lower threshold
$L\in\{0,1,\ldots,2^k-1\}$, define the $k$-bit lower-bound penalty $H_k(X,L) \equiv H_{\mathrm{lower}}(X,L)$ recursively as follows. For $k=1$,
\[
    H_1(X,0)=0,
    \qquad
    H_1(X,1)=1-x_1.
\]
For the extension from $k$ to $k+1$ bits, write
\[
    X_{k+1}=2^k x_{k+1}+X_k
\]
and
\[
    L=2^k\ell_{k+1}+L',
\]
where $x_{k+1},\ell_{k+1}\in\{0,1\}$ and
$0\leq X_k,L'<2^k$, and define
\[
H_{k+1}(X_{k+1},L)
=
\begin{cases}
(1-x_{k+1})H_k(X_k,L'),
& \ell_{k+1}=0,\\[4pt]
(1-x_{k+1})+H_k(X_k,L'),
& \ell_{k+1}=1.
\end{cases}
\]
Then $H_k(X,L)$ is a nonnegative integer-valued polynomial in the
original binary variables satisfying
\[
    H_k(X,L)=0
    \quad\Longleftrightarrow\quad
    L\leq X.
\]
\end{theorem}

\begin{proof}
We prove the stated properties by mathematical induction on the
bit width $k$.

\textit{Base case.}
Let $k=1$. Then $X=x_1$, and the possible lower thresholds are
$L=0$ and $L=1$.

For $L=0$, the recursive definition gives
\[
    H_1(X,0)=0.
\]
Both representable values satisfy $0\leq X$, so the required
property holds.

For $L=1$, the recursive definition gives
\[
    H_1(X,1)=1-x_1.
\]
Hence $H_1(X,1)=0$ exactly when $X=1$, and it is strictly positive
exactly when $X=0$. Therefore, the theorem holds for $k=1$.

\textit{Inductive hypothesis.}
Assume that for some $k\geq 1$ and for every lower threshold
$L'\in\{0,1,\ldots,2^k-1\}$, the recursively constructed
polynomial $H_k(X_k,L')$ is nonnegative and integer-valued and
satisfies
\[
    H_k(X_k,L')=0
    \quad\Longleftrightarrow\quad
    L'\leq X_k,
\]
and
\[
    H_k(X_k,L')>0
    \quad\Longleftrightarrow\quad
    X_k<L'.
\]

\textit{Inductive step.}
Consider an arbitrary $(k+1)$-bit integer and an arbitrary
$(k+1)$-bit lower threshold written as
\[
    X_{k+1}=2^kx_{k+1}+X_k
\]
and
\[
    L=2^k\ell_{k+1}+L',
\]
where $x_{k+1},\ell_{k+1}\in\{0,1\}$ and
$0\leq X_k,L'<2^k$.

\textit{Case A: $\ell_{k+1}=0$.}
Then $L=L'$. By the recursive definition,
\begin{equation}
    H_{k+1}(X_{k+1},L)
    =
    (1-x_{k+1})H_k(X_k,L').
    \label{eq:lower_caseA}
\end{equation}

If $x_{k+1}=1$, then
\[
    X_{k+1}\geq 2^k>L,
\]
and simultaneously
\[
    H_{k+1}(X_{k+1},L)=0.
\]

If $x_{k+1}=0$, then
\[
    H_{k+1}(X_{k+1},L)
    =
    H_k(X_k,L'),
\]
while
\[
    L\leq X_{k+1}
    \quad\Longleftrightarrow\quad
    L'\leq X_k.
\]
By the inductive hypothesis,
\[
    H_{k+1}(X_{k+1},L)=0
    \quad\Longleftrightarrow\quad
    L\leq X_{k+1}.
\]
Likewise,
\[
    H_{k+1}(X_{k+1},L)>0
    \quad\Longleftrightarrow\quad
    X_{k+1}<L.
\]
Thus the required property holds in Case A.

\textit{Case B: $\ell_{k+1}=1$.}
Then $L=2^k+L'$. By the recursive definition,
\begin{equation}
    H_{k+1}(X_{k+1},L)
    =
    (1-x_{k+1})+H_k(X_k,L').
    \label{eq:lower_caseB}
\end{equation}
Both terms are nonnegative. Therefore,
\[
    H_{k+1}(X_{k+1},L)=0
\]
if and only if
\[
    x_{k+1}=1
    \quad\text{and}\quad
    H_k(X_k,L')=0.
\]
By the inductive hypothesis, this is equivalent to
\[
    x_{k+1}=1
    \quad\text{and}\quad
    L'\leq X_k.
\]
Since $L=2^k+L'$, these conditions are equivalent to
\[
    L\leq X_{k+1}.
\]
If either condition fails, at least one of the two nonnegative
terms is strictly positive, and hence
\[
    H_{k+1}(X_{k+1},L)>0.
\]
Thus,
\[
    H_{k+1}(X_{k+1},L)=0
    \quad\Longleftrightarrow\quad
    L\leq X_{k+1}
\]
in Case B.

In both cases, $H_{k+1}$ is obtained from the nonnegative
integer-valued polynomial $H_k$ by addition or multiplication with
the binary literal $1-x_{k+1}$. Hence it remains a nonnegative
integer-valued polynomial, and its zero set is exactly
$\{X_{k+1}:L\leq X_{k+1}\}$.

Therefore, the stated properties hold for bit width $k+1$ whenever
they hold for bit width $k$. By mathematical induction, the theorem
holds for every integer $k\geq 1$.
\end{proof}

\begin{corollary}[Explicit finite sum-product representation]
\label{cor:lower_explicit}
Let
\[
    L=(\ell_k\ell_{k-1}\cdots\ell_1)_2.
\]
The recursively constructed lower-bound penalty in
Theorem~\ref{thm:lower} admits the equivalent explicit
finite sum-product representation
\begin{equation}
    \Hlower(X,L)
    =
    \sum_{\substack{1\leq i\leq k\\\ell_i=1}}
    (1-x_i)
    \prod_{\substack{j>i\\\ell_j=0}}(1-x_j).
    \label{eq:lower_explicit}
\end{equation}
An empty product is defined as one, and an empty sum is defined as
zero. Hence, each one bit $\ell_i=1$ of the threshold corresponds
to exactly one product term.
\end{corollary}

\begin{proof}
Starting from the one-bit construction and successively appending higher-order threshold bits, a one threshold bit creates a new product term initialized by the corresponding complemented literal, whereas a zero threshold bit multiplies every existing product term by the corresponding complemented literal. Repeated application of these two rules gives Eq.~\eqref{eq:lower_explicit}.
\end{proof}

As a concrete illustration of the base case and one inductive step,
Table~\ref{tab:lower_k2} lists the resulting penalty functions for all
two-bit lower thresholds.

\begin{table}[H]
\centering
\caption{Lower-bound penalty functions for all thresholds at $k=2$.}
\label{tab:lower_k2}
\begin{tabular}{cc}
\toprule
Lower threshold $L$ & $\Hlower(X,L)$ \\
\midrule
$0=(00)_2$ & $0$ \\
$1=(01)_2$ & $(1-x_2)(1-x_1)$ \\
$2=(10)_2$ & $(1-x_2)$ \\
$3=(11)_2$ & $(1-x_2)+(1-x_1)$ \\
\bottomrule
\end{tabular}
\end{table}

Each function in Table~\ref{tab:lower_k2} vanishes exactly for the
assignments satisfying $L\leq X$ and is strictly positive for those
satisfying $X<L$.

Table~\ref{tab:lower_k3} extends the construction to three bits. Each
entry is obtained from the corresponding two-bit function by
applying Case A or Case B of the inductive step.

\begin{table}[H]
\centering
\caption{Lower-bound penalty functions for all thresholds at $k=3$.}
\label{tab:lower_k3}
\begin{tabular}{ccc}
\toprule
$(\ell_3\ell_2\ell_1)_2$ & $L$ & $\Hlower(X,L)$ \\
\midrule
000 & 0 & $0$ \\
001 & 1 & $(1-x_3)(1-x_2)(1-x_1)$ \\
010 & 2 & $(1-x_3)(1-x_2)$ \\
011 & 3 & $(1-x_3)\bigl((1-x_2)+(1-x_1)\bigr)$ \\
100 & 4 & $(1-x_3)$ \\
101 & 5 & $(1-x_3)+(1-x_2)(1-x_1)$ \\
110 & 6 & $(1-x_3)+(1-x_2)$ \\
111 & 7 & $(1-x_3)+(1-x_2)+(1-x_1)$ \\
\bottomrule
\end{tabular}
\end{table}

The table also illustrates the reuse of lower-bit penalty structures.

\paragraph*{Algorithm~2: Iterative construction of
  \texorpdfstring{$\Hlower$}{H\_lower}}\mbox{}\\
Given the binary representation
\[
    L=(\ell_k\ell_{k-1}\cdots\ell_2\ell_1)_2,
\]
the lower-bound penalty function $\Hlower(X,L)$ is constructed as
follows.
\begin{enumerate}
    \item If $\ell_1=0$, initialize $H = 0$.
          If $\ell_1=1$, initialize $H = 1-x_1$.

    \item For $i=2,3,\ldots,k$, update
    \[
        H\leftarrow
        \begin{cases}
            (1-x_i)H, & \text{if } \ell_i=0,\\
            (1-x_i)+H, & \text{if } \ell_i=1.
        \end{cases}
    \]
\end{enumerate}
The final value of $H$ is $\Hlower(X,L)$. Algorithm~2 is the
iterative form of the recurrence proved in Theorem~\ref{thm:lower},
so its correctness follows from that theorem.

\paragraph*{Example $(k=6,\ L=17)$}\mbox{}\\
For
\[
    L=17=(010001)_2,
\]
Algorithm~2 gives
\begin{equation}
    \Hlower(X,17)
    =
    (1-x_6)
    \bigl(
        (1-x_5)
        +(1-x_4)(1-x_3)(1-x_2)(1-x_1)
    \bigr).
    \label{eq:example_lower}
\end{equation}
This function has product-term count two and degree five. In contrast,
the trivial construction contains seventeen product terms, one for
each infeasible integer in $\{0,1,\ldots,16\}$, and has degree six.
Thus, the proposed construction reduces the product-term count from
seventeen to two and the degree from six to five.

\subsection{Structural Duality}
\label{subsec:duality}

The upper-bound and lower-bound penalty functions are related by
one's-complement duality. Complementing every bit of $X$ gives
\[
    \overline{X}
    =
    \sum_{i=1}^{k}(1-x_i)2^{i-1}
    =
    2^k-1-X.
\]
Similarly, complementing every bit of the lower threshold $L$ gives
\[
    \overline{L}=2^k-1-L.
\]
The lower-bound condition is then transformed into the upper-bound
condition
\[
    L\leq X
    \quad\Longleftrightarrow\quad
    \overline{X}\leq\overline{L}.
\]

Accordingly, the two penalty functions satisfy
\begin{equation}
    \Hlower(X,L)
    =
    \Hupper(2^k-1-X,\,2^k-1-L).
    \label{eq:structural_duality}
\end{equation}
Under this substitution, each literal $x_i$ in the upper-bound
construction is replaced by $1-x_i$, while each zero bit of the
complemented threshold corresponds to a one bit of $L$.

This relation is a substitution on the encoded integer and threshold, not a Boolean negation of the penalty function such as $1-\Hupper$.

Although the lower-bound construction can also be derived from this duality,
Theorem~\ref{thm:lower} provides an independent inductive verification of the lower-bound recurrence.

\subsection{Degree and Product-Term Count}
\label{subsec:degree_count}

We now compare the product-term count and degree of the proposed and
trivial constructions. As defined in
Section~\ref{subsec:prelim}, both quantities are evaluated on the
expressions generated by the respective constructions,
before algebraic simplification. The endpoint cases $U=2^k-1$ and
$L=0$ produce the zero function and are omitted only from the degree
comparisons.

For the upper-bound penalty function, the explicit representation in
Eq.~\eqref{eq:upper_explicit} contains one product term for each zero
bit of $U$. Its product-term count is therefore at most $k$. The
trivial construction instead contains one product term for every
integer that violates $X\leq U$, giving $2^k-U-1$ product terms. The
proposed construction has strictly fewer product terms for every
$U$ from $0$ to $2^k-2$, except at $U=2^k-2$, where both
constructions contain one product term.

The degree of the upper-bound construction is determined by the least
significant zero bit of $U$ and the one bits above it. At
$U=2^k-2$, the proposed penalty function is
$x_kx_{k-1}\cdots x_1$. It therefore has degree $k$ and coincides
with the trivial construction. For every smaller threshold, the
proposed construction has degree less than $k$, whereas every product
term generated by the trivial construction has degree $k$.

For the lower-bound penalty function, the explicit representation in
Eq.~\eqref{eq:lower_explicit} contains one product term for each one
bit of $L$. Its product-term count is therefore also at most $k$.
The trivial construction contains one product term for every integer
from $0$ to $L-1$, giving $L$ product terms. The proposed
construction has strictly fewer product terms for every positive
threshold except $L=1$, where both constructions contain one product
term.

The degree of the lower-bound construction is determined by the least
significant one bit of $L$ and the zero bits above it. At $L=1$, the
proposed penalty function is
$(1-x_k)(1-x_{k-1})\cdots(1-x_1)$. It therefore has degree $k$ and
coincides with the trivial construction. For every threshold greater
than one, the proposed construction has degree less than $k$, whereas
every product term generated by the trivial construction has degree
$k$.

In terms of the representable range size $R$ introduced in Section~\ref{sec:intro}, the trivial construction
can require $O(R)$ product terms, whereas the
proposed construction requires at most $\log_2 R$ product terms and
therefore scales as $O(\log R)$. Averaged uniformly over all representable thresholds, the
product-term count is $\frac{1}{2}\log_2 R$ for the proposed
construction, compared with $\frac{R-1}{2}$ for the trivial
construction. Thus, the average product-term count scales logarithmically rather than linearly with $R$.

\begin{figure}[H]
    \centering
    \includegraphics[width=0.88\linewidth]
        {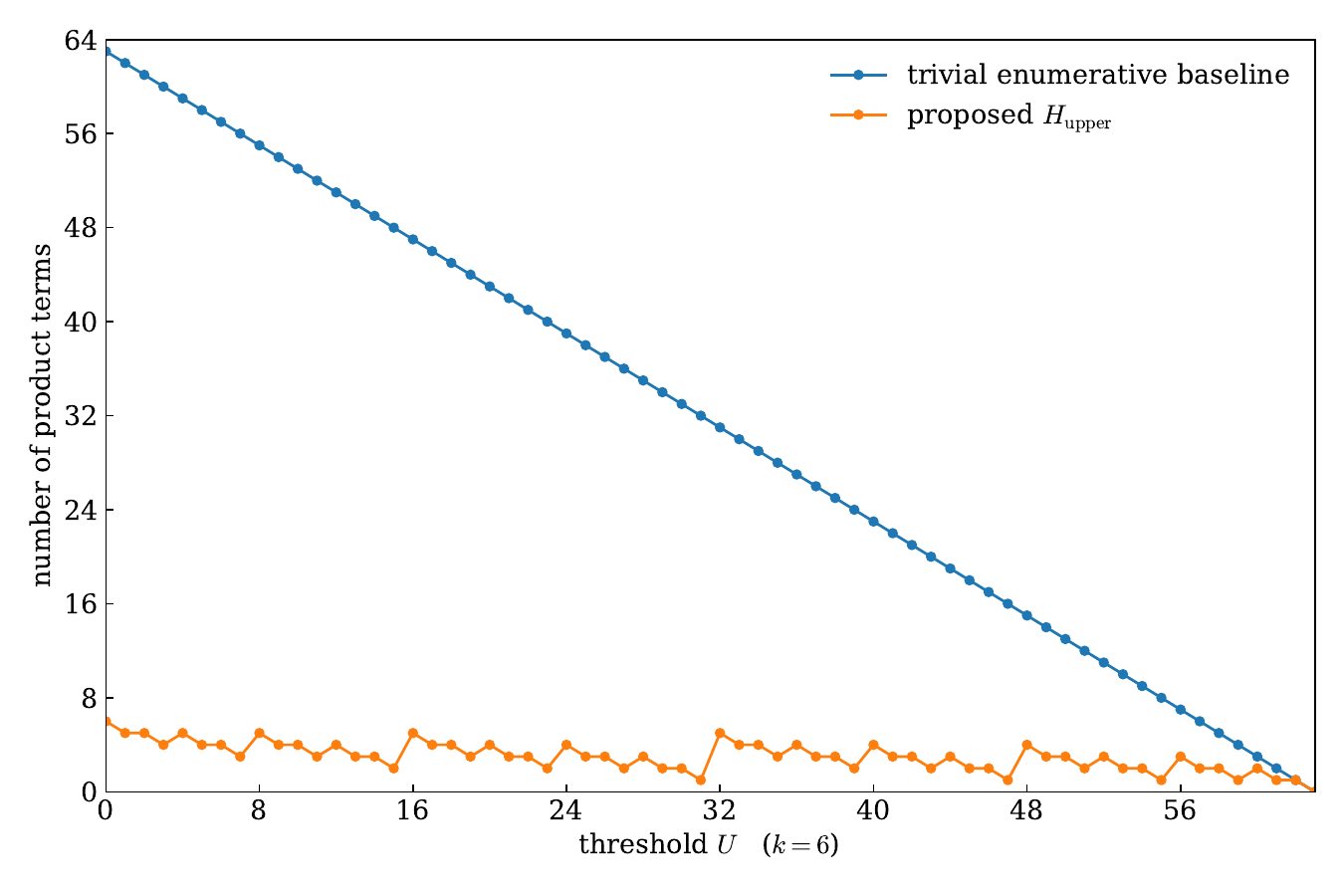}
    \caption{Product-term counts of the proposed and trivial
    upper-bound constructions for $k=6$. The proposed count equals
    the number of zero bits in $U$, whereas the trivial count is
    $2^6-U-1$.}
    \label{fig:product_term_count}
\end{figure}

Figure~\ref{fig:product_term_count} illustrates these product-term counts for the upper-bound constructions at $k=6$. The proposed count remains at most six, whereas the trivial count ranges up to $63$; the two coincide at $U=62$ and at the vacuous endpoint $U=63$.

\begin{figure}[H]
    \centering
    \includegraphics[width=0.88\linewidth]
        {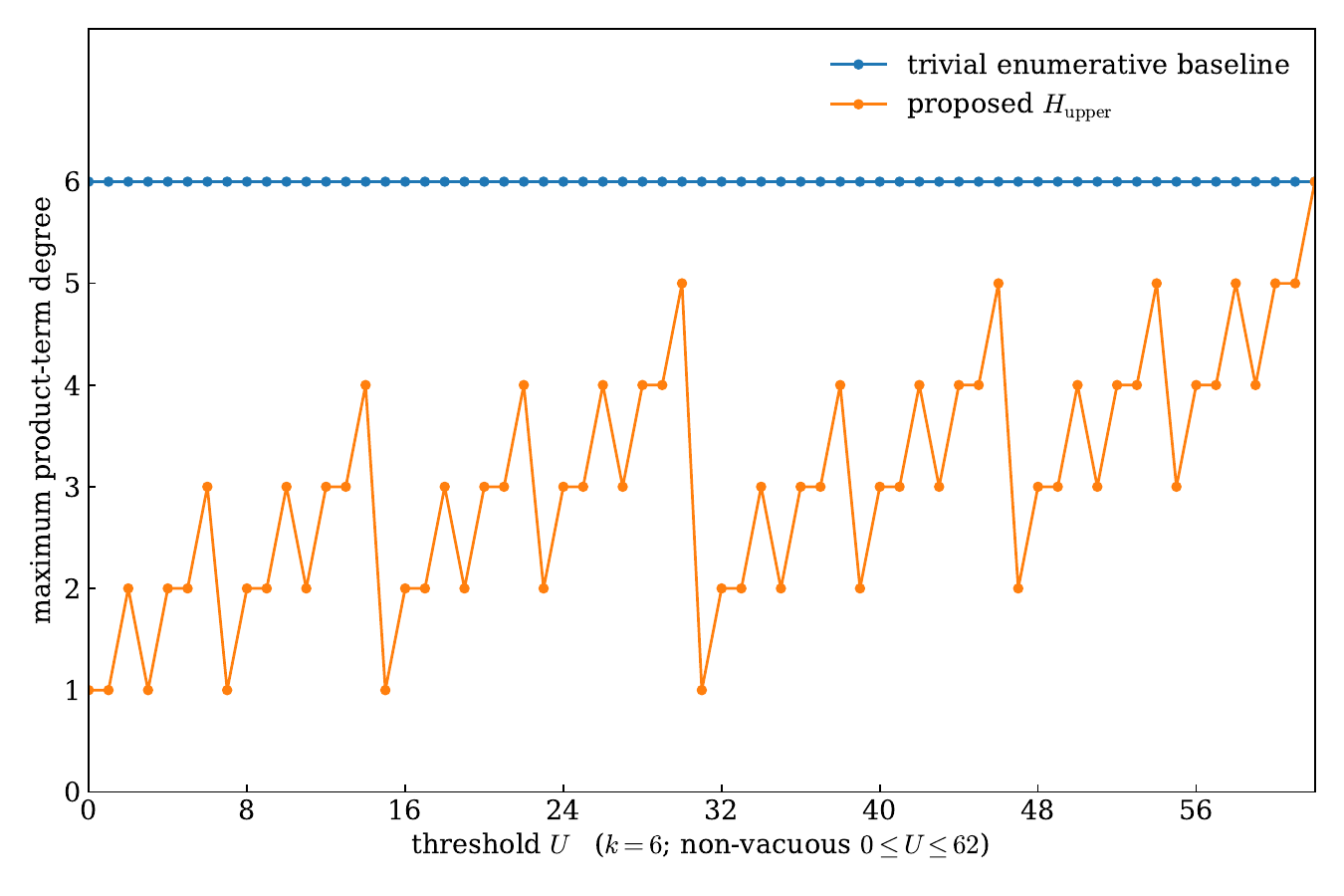}
    \caption{Degrees of the proposed and trivial upper-bound
    constructions for $U=0,\ldots,62$ at $k=6$. The trivial
    construction has degree six throughout, whereas the proposed
    construction reaches degree six only at $U=62$.}
    \label{fig:degree_comparison}
\end{figure}

Figure~\ref{fig:degree_comparison} illustrates the corresponding degree comparison at $k=6$. The proposed degree remains below six for
$U=0,\ldots,61$ and reaches six only at $U=62$, where the two constructions coincide.

\section{Penalty-Function Arithmetic Composition}
\label{sec:logic}

\begin{theorem}[Interval penalty function]
\label{thm:interval}
Let $0\leq L\leq U<2^k$. The penalty function for the interval
constraint $L\leq X\leq U$ is
\begin{equation}
    \Hin(X,L,U)
    =
    \Hlower(X,L)+\Hupper(X,U).
    \label{eq:Hin}
\end{equation}
\end{theorem}

\begin{proof}
The interval constraint is satisfied exactly when both
\[
    L\leq X
    \quad\text{and}\quad
    X\leq U
\]
hold. We therefore require the composite penalty function to be zero
exactly when both the lower bound and upper bound penalty functions
are zero. Since both functions are nonnegative, addition has precisely
this property. Hence the stated sum has exactly the required zero set.
\end{proof}

\begin{theorem}[Exclusion penalty function]
\label{thm:exclusion}
Let $0\leq L\leq U<2^k$. The penalty function for the exclusion
constraint
\[
    X\leq L
    \quad\text{or}\quad
    U\leq X
\]
is
\begin{equation}
    \Hout(X,L,U)
    =
    \Hupper(X,L)\Hlower(X,U).
    \label{eq:Hout}
\end{equation}
\end{theorem}

\begin{proof}
The exclusion constraint is satisfied exactly when
\[
    X\leq L
    \quad\text{or}\quad
    U\leq X.
\]
By Theorems~\ref{thm:upper} and~\ref{thm:lower}, these two conditions
are equivalent to
\[
    \Hupper(X,L)=0
    \quad\text{or}\quad
    \Hlower(X,U)=0.
\]
The composite penalty must therefore be zero whenever at least one
component penalty is zero and must be positive when both component
penalties are positive. Multiplication provides this zero set.

If $X\leq L$, then $\Hupper(X,L)=0$, and hence
\[
    \Hout(X,L,U)=0.
\]
Similarly, if $U\leq X$, then $\Hlower(X,U)=0$, and again
\[
    \Hout(X,L,U)=0.
\]

Conversely, suppose that
\[
    \Hout(X,L,U)=0.
\]
Since
\[
    \Hout(X,L,U)
    =
    \Hupper(X,L)\Hlower(X,U),
\]
at least one factor must be zero. Theorems~\ref{thm:upper}
and~\ref{thm:lower} then imply
\[
    X\leq L
    \quad\text{or}\quad
    U\leq X.
\]

Finally, if $L<X<U$, then
\[
    \Hupper(X,L)>0
    \quad\text{and}\quad
    \Hlower(X,U)>0.
\]
Their product is therefore strictly positive. Hence
$\Hout(X,L,U)$ is zero exactly on the required feasible set and is
strictly positive throughout the excluded region.
\end{proof}

Thus, for nonnegative penalty functions, addition realizes the
intersection of component zero sets, whereas multiplication realizes
their union.

The arithmetic composition also determines the degree of the resulting
penalty function. Because the interval penalty function is formed by
addition, its degree is no greater than the larger degree of its two
component functions. Because the exclusion penalty function is formed
by multiplication, its degree is no greater than the sum of their
degrees.

\begin{figure}[H]
    \centering
    \includegraphics[width=\columnwidth]
    {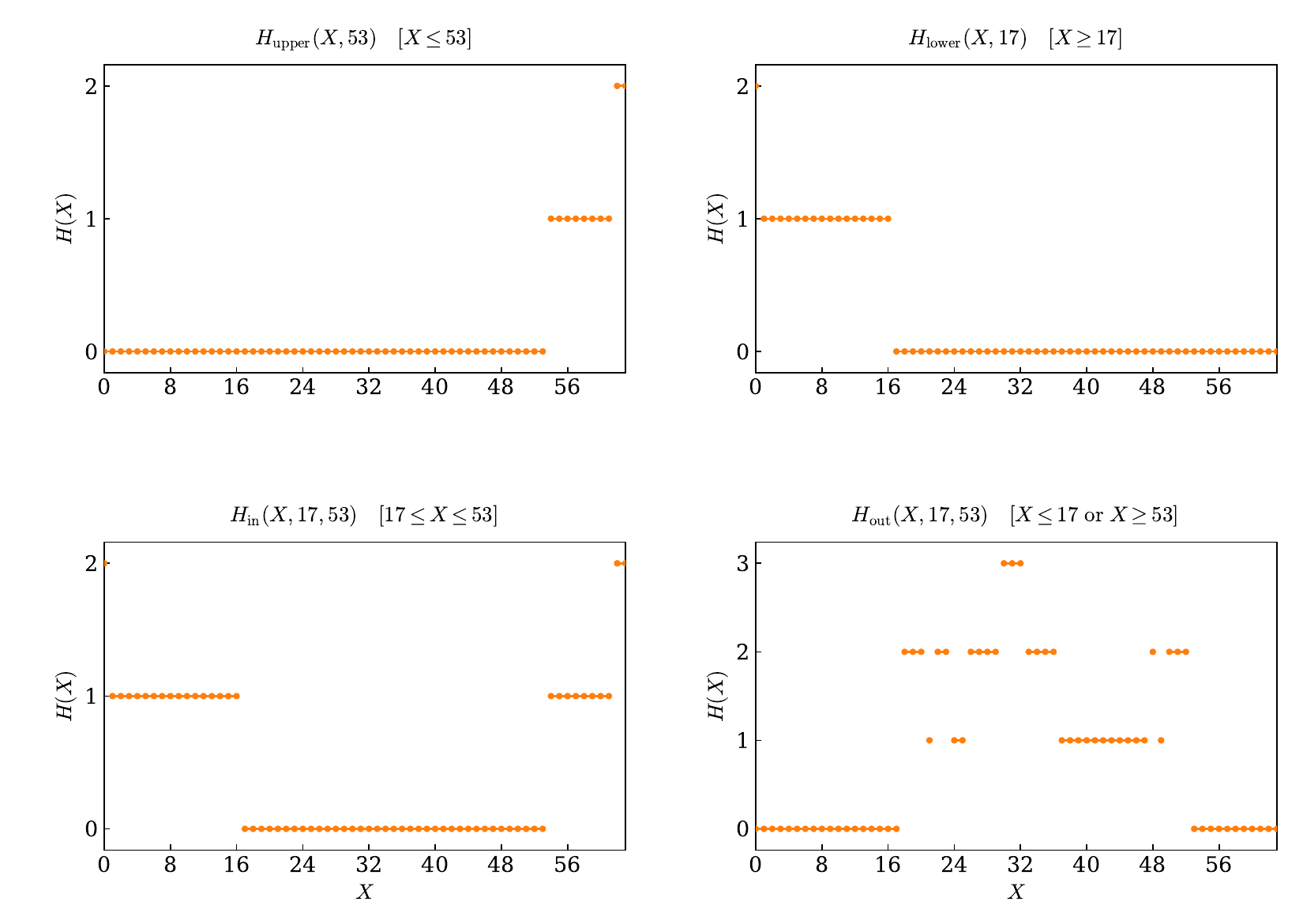}
    \caption{The four canonical penalty functions evaluated over
    $X\in\{0,\ldots,63\}$ for $k=6$, $L=17$, and $U=53$.
    The upper-bound and lower-bound functions are shown in the top
    row. Their pointwise sum gives the interval penalty function in
    the bottom-left panel, which is zero exactly for
    $17\leq X\leq53$. The pointwise product
    $\Hupper(X,17)\Hlower(X,53)$ gives the exclusion penalty function
    in the bottom-right panel, which is zero exactly for $X\leq17$
    or $53\leq X$. Every function is strictly positive outside its
    specified feasible set.}
    \label{fig:landscape}
\end{figure}

Figure~\ref{fig:landscape} illustrates the four canonical penalty functions for $k=6$, $L=17$, and $U=53$. The top row shows the two single-bound penalties. Their pointwise sum yields the interval penalty, whose zero set is $17\leq X \leq 53$, while the product
$\Hupper(X,17)\Hlower(X,53)$ yields the exclusion penalty, whose zero set is $X\leq17$ or
$53\leq X$.

\section{Discussion}
\label{sec:discussion}

\subsection{Numerical Verification}
\label{subsec:numerical_verification}

As a numerical consistency check, we exhaustively evaluated all 64 assignments for $k=6$ using $L=17$ and $U=53$. All four penalty functions reproduced the zero sets predicted by Theorems~\ref{thm:upper}--\ref{thm:exclusion}, with strictly positive values elsewhere.

The Python code provided as Supplementary Material reproduces this exhaustive evaluation for the upper-bound, lower-bound, interval, and exclusion penalty functions.

\subsection{Mathematical Generality and Implementation Implications}
\label{subsec:implications}

The proposed penalty functions address the four structural limitations of the trivial construction identified in Section~\ref{subsec:motivation}. We discuss these improvements in the
same order: growth in product-term count, fixed maximum degree, lack
of hierarchical structure, and poor mathematical reusability and
extensibility. We also consider their implications for classical,
quantum-inspired, and quantum optimization implementations.

\textbf{(i) Reduction in product-term count.}\mbox{}\\
The hierarchical construction reuses the penalty function obtained for the lower bit positions instead of introducing one product term for each infeasible integer. As shown in Section~\ref{subsec:degree_count}, the proposed single-bound penalties contain at most $k=O(\log R)$ product terms, whereas the trivial construction can require $O(R)$ terms. This reduction complements the compact variable representation already provided by binary encoding. This may reduce the number of terms that must be stored or evaluated in classical and quantum-inspired solvers, and may also reduce the number of higher-order interaction terms that must be represented in quantum implementations.

\textbf{(ii) Reduction in maximum degree.}\mbox{}\\
As established in Section~\ref{subsec:degree_count}, the proposed single-bound penalties have lower maximum degree than the trivial construction over nearly the entire non-vacuous threshold range. A
lower degree reduces the number of binary variables that occur
simultaneously in a single product term. This may decrease the evaluation burden in classical and quantum-inspired solvers that operate on higher-order binary polynomials, including HOBOTAN~\cite{HOBOTAN2024}. It may also simplify the representation of these penalty functions on future quantum hardware supporting higher-order interactions. On current quantum annealers, where higher-order terms generally require quadratization, the practical benefit depends on the selected reduction method and hardware connectivity~\cite{Mandal2020,Dattani2019,Konz2021,Abbas2024}. 

\textbf{(iii) Reuse of hierarchical structure.}\mbox{}\\
 The proposed construction forms each $(k+1)$-bit penalty from the corresponding $k$-bit penalty by one additional additive or multiplicative step, thereby preserving and reusing the lower-bit structure. In contrast, the trivial construction treats each infeasible bit string independently and shares no substructure among the resulting product terms. This recursive organization is structurally analogous to hierarchical constructions used elsewhere in quantum computing~\cite{Barenco1995}.

\textbf{(iv) Mathematical reusability and extensibility.}\mbox{}\\
Beyond the hierarchical reuse described above, the proposed construction provides a common algebraic framework for combining threshold penalties. Because the single-bound penalties are nonnegative, addition produces interval penalties and multiplication produces exclusion penalties through their zero-set behavior. These constructions can therefore be reused and combined without returning to term-by-term enumeration, providing a basis for extensions to more complex inequality constraints.

Together, these properties provide a more compact and reusable representation of binary inequality penalties than the trivial construction.

\subsection{Additional Advantages}
\label{subsec:additional}

Beyond the reductions in product-term count and degree, the proposed
penalty functions offer three additional advantages.

First, the proposed upper-bound penalty provides a compact way to exclude
unused bit strings when binary encoding is applied to a domain whose number
of admissible values is not a power of two. If a variable has $N$ admissible
values, a binary representation requires only $\lceil \log_2 N\rceil$ bits.
When $N$ is not a power of two, however, some bit strings represent integers
greater than $N-1$ and therefore do not correspond to valid states of the
original variable. In a higher-order formulation, these unused bit strings
can be excluded by a trivial penalty function that enumerates them
individually, but such a construction inherits the structural limitations
identified in Section~\ref{subsec:motivation}. The proposed upper-bound
penalty function $\Hupper(X,N-1)$ excludes the same unused values through
the hierarchical threshold-based construction developed in this paper.
Similar unused bit strings can arise when compact binary encodings are used
for categorical choices in problems such as the Travelling Salesman Problem
(TSP). The method may therefore be useful not only for constraints originally
stated as inequalities, but also for enforcing the valid domain of compact
binary encodings.

Second, the single-bound penalty functions have a bounded range of
penalty values. Their smallest positive value is one, and their
largest value is at most the bit width $k$, because each product term
takes a binary value and each construction contains at most $k$
product terms. This bounded range makes the scale of the single-bound
penalties predictable and may simplify the relative weighting between
the penalty term, other constraints, and the objective function. It
does not, however, eliminate the need to choose an appropriate overall
penalty coefficient.

Third, the same construction can be applied to fractional-valued quantities by interpreting the binary-encoded integer at a fixed numerical scale. The inequality thresholds are represented at the same scale, allowing the upper-bound, lower-bound, interval, and exclusion penalty functions to be used without changing their underlying recursive construction. This makes the framework applicable not only to integer-valued variables but also to variables that require fractional resolution. Greater numerical precision can be obtained by increasing the bit width while retaining the same recursive rules.

\subsection{Limitations and Outlook}
\label{subsec:outlook}

The first limitation is that the degree reduction is not uniform in the bit width. As shown in Section~\ref{subsec:degree_count}, the proposed and trivial constructions coincide in degree at one non-vacuous boundary threshold for each bound. More importantly, the degree of the proposed construction is not bounded by a constant independent of $k$ and can still increase as the bit width $k$ increases.

The second limitation concerns implementation on current quantum
annealers. Although the proposed penalty functions introduce no
auxiliary variables in their original higher-order form, they can
contain product terms of degree greater than two. On hardware restricted
to two-body interactions, such terms generally require quadratization,
which can reintroduce auxiliary variables and embedding overhead
depending on the reduction method and hardware connectivity~\cite{Dattani2019,Mandal2020}. The slack-free property therefore
applies to the original higher-order formulation but may not be fully
preserved after quadratization. Future hardware supporting native
multi-body interactions could mitigate this limitation by reducing or
eliminating the need for quadratization.

The third limitation concerns the penalty values assigned to infeasible
assignments. Although the proposed penalty functions distinguish feasible
and infeasible assignments exactly, their positive values do not generally
increase with the magnitude of the constraint violation. A larger violation
therefore need not receive a larger penalty. Consequently, the penalty
values do not provide a monotonic measure of violation severity or
directional guidance toward the feasible region. The construction is
therefore not intended as a violation-distance penalty for optimization
methods that rely on such guidance.

The fourth limitation concerns empirical optimization performance. The
present work establishes the correctness of the proposed penalty functions,
analyzes their product-term counts and degrees, and verifies their zero sets
exhaustively for the representative $k=6$ case. It does not, however,
compare their optimization performance with slack-variable encodings,
unbalanced penalization~\cite{MontanezBarrera2024}, or the trivial
construction on application-level instances. The reductions in product-term
count and degree should therefore not be interpreted as guarantees of a
higher feasible-solution rate, shorter time to solution, or lower
implementation cost after quadratization. Systematic comparisons using
common solvers, quadratization methods, penalty coefficients, and problem
instances remain an important direction for future work.

As discussed in Section~\ref{subsec:additional}, the same
hierarchical construction can be applied to quantities represented at a
fixed fractional scale. This property opens a path toward applications
that require substantially finer numerical resolution, including
optimization problems in chemistry and other scientific fields. By
increasing the bit width $k$, the upper-bound, lower-bound, interval,
and exclusion penalty functions can be reused at progressively finer
scales without redesigning the underlying recursive rules. Future work
may therefore extend the present framework from integer-valued
constraints to higher-precision scientific optimization models and
investigate how the required resolution, polynomial degree, and solver
performance interact in practical applications.

A further direction is to extend the arithmetic composition developed in
Section~\ref{sec:logic} beyond the interval and exclusion constraints
considered here. Addition and multiplication may provide a basis for
combining larger collections of nonnegative penalty functions. Extensions
to complement operations, exclusive alternatives such as XOR, and other
compound constraint structures could also be investigated. Such
developments may broaden the applicability of the framework to more complex
binary-encoded constrained optimization problems.

\section{Conclusions}
\label{sec:conclusion}

This paper has presented exact and slack-free penalty functions for
upper-bound, lower-bound, interval, and exclusion constraints on
binary-encoded integers. The upper-bound and lower-bound penalties are constructed hierarchically from the threshold bits, and equivalent explicit finite sum-product representations are obtained by expanding the underlying recursive rules.

The resulting inequality penalty functions can be generated from the binary digits of their thresholds without enumerating
infeasible assignments or introducing slack variables in the original
higher-order formulation. Relative to the trivial construction, the proposed single-bound penalties reduce both product-term count and degree at all but one non-vacuous boundary threshold for each bound, while both constructions reduce to the zero polynomial at the vacuous endpoints. The arithmetic composition of the upper-bound and lower-bound functions
further provides a basis for constructing compound constraints such as
intervals and exclusions, and may support extensions to more complex
logical compositions.
These results provide a mathematically reusable framework for exact
inequality penalties and may support broader constrained formulations
in higher-order binary optimization.

\section*{CRediT authorship contribution statement}

\textbf{Naoaki Mochida:} Conceptualization; Methodology; Formal
analysis; Software; Validation; Visualization; Writing -- original
draft; Writing -- review \& editing.

\textbf{Shoya Yasuda:} Conceptualization; Methodology; Validation;
Supervision; Writing -- review \& editing.

\section*{Funding}

This research did not receive any specific grant from funding agencies in the public, commercial, or not-for-profit sectors.

\section*{Declaration of competing interest}

The authors declare no competing
interests or personal relationships that could have appeared to
influence the work reported in this paper.

\section*{Data availability}

No external datasets were used or generated in this theoretical
study. The Python code used for the numerical verification
described in Section~\ref{subsec:numerical_verification} is
provided as Supplementary Material.

\section*{Declaration of generative AI and AI-assisted technologies in the manuscript preparation process}

During the preparation of this work, the authors used ChatGPT (OpenAI) to improve language and readability. After using this tool, the authors reviewed and edited the content as needed and take full responsibility for the content of the publication. All scientific claims, mathematical derivations, results, and conclusions were developed and verified by the authors.

\bibliographystyle{elsarticle-num}
\bibliography{refs}
\end{document}